\documentclass[a4paper,12pt,reqno]{amsart}
\usepackage{graphics}
\newtheorem{lem}{Lemma}

\newtheorem{theo}{Theorem}

\newtheorem{cor}{Corollary}

\numberwithin{equation}{section}

\newcommand{\sgn}{\operatorname{sgn}}

\newcommand{\id}{\operatorname{id}}

\newcommand{\Op}{\operatorname{Op}}

\begin{document}

\title[Monotone triangles]{The number of monotone triangles with prescribed bottom row}

\author[Ilse Fischer]{\box\Adr}

\newbox\Adr
\setbox\Adr\vbox{ \centerline{ \large Ilse Fischer} \vspace{0.3cm}
\centerline{Fakult\"at f\"ur Mathematik, Universit\"at Wien}
\centerline{Nordbergstrasse 15, A-1090 Wien, Austria}
\centerline{E-mail: {\tt Ilse.Fischer@univie.ac.at}} }

\begin{abstract} We show that the number of monotone triangles with prescribed bottom row 
$(k_1,\ldots,k_n) \in \mathbb{Z}^n$, $k_1 < k_2 < \ldots < k_n$, is given by a simple product 
formula which remarkably involves (shift) operators. Monotone triangles with bottom row $(1,2,\ldots,n)$ 
are in bijection with $n \times n$ alternating sign matrices.   
\end{abstract}

\maketitle

\section{Introduction}

An {\it alternating sign matrix} is a square matrix of $0$s, $1$s and $-1$s for which 
the sum of entries in each row and in each column is $1$ and the non-zero entries of each row and 
of each column alternate in sign. For instance, 
$$
\left(
\begin{array}{rrrrr}
0 & 0 & 0 & 1 & 0 \\
0 & 1 & 0 & -1& 1 \\
1 & -1& 1 & 0 & 0 \\
0 & 1 & -1 & 1 & 0 \\
0 & 0 &  1 & 0 & 0
\end{array}
\right)
$$
is an alternating sign matrix. In the early 1980s, Robbins and Rumsey \cite{robbins} introduced
alternating sign matrices in the course of generalizing a determinant evaluation algorithm. 
Out of curiosity they posed the question for the number of alternating sign matrices 
of fixed size and, together with Mills, they came up with the appealing conjecture \cite{mills} that the number of 
$n \times n$ alternating sign matrices is 
\begin{equation}
\label{asmformula}
\prod_{j=1}^{n} \frac{(3j-2)!}{(n+j-1)!}.
\end{equation}
This turned out to be one of the hardest problems in enumerative combinatorics within the last decades. 
In 1996 Zeilberger \cite{zeilberger} finally succeeded in proving their conjecture. Then, some 
months later, Kuperberg \cite{kuperberg} realized that alternating sign matrices are equivalent to a 
model in statistical physics for two-dimensional square ice. Using a determinantal expression for the 
partition function of this model discovered earlier by physicists, he was able to provide a shorter proof of the formula. 
For a nice exposition on this topic see~\cite{bressoud}.

Alternating sign matrices can be translated into certain triangular arrays of positive integers, 
called {\it monotone triangles}. Monotone triangles are probably the right guise 
of alternating sign matrices for a recursive treatment \cite[Section 2.3]{bressoud}.
In order to obtain the monotone triangle corresponding to a given alternating sign matrix,
replace every entry in the matrix by the sum of elements in the same column above, the 
entry itself included. In our running example we obtain
$$
\left(
\begin{array}{rrrrr}
0 & 0 & 0 & 1 & 0 \\
0 & 1 & 0 & 0& 1 \\
1 & 0& 1 & 0 & 1 \\
1 & 1 & 0 & 1 & 1 \\
1 & 1 &  1 & 1 & 1
\end{array}
\right).
$$
Row by row we record the columns that contain a $1$ and obtain the following triangular array.
$$
\begin{array}{ccccccccc}
  &  &  & & 4  &  &  &  &  \\
  &  &  &2 &   &  5&  &  &  \\
  &  &1  & &  3 &  & 5 &  &  \\
  &1  &  & 2&   & 4 &  &  5&  \\
1 & & 2 & & 3 & & 4 & & 5
\end{array}
$$
This is the monotone triangle corresponding to the alternating sign matrix above. Observe that 
it is weakly increasing in northeast direction and in southeast direction. Moreover, it is strictly 
increasing along rows. In general, a monotone triangle with $n$ rows is a triangular array 
$(a_{i,j})_{1 \le j \le i \le n}$ of integers such that $a_{i,j} \le a_{i-1,j} \le a_{i,j+1}$ 
and $a_{i,j} < a_{i,j+1}$ for all $i, j$. It is not too hard to see that 
monotone triangles with $n$ rows and bottom row $(1, 2, \ldots, n)$, i.e. $a_{n,j}=j$, 
are in bijection with $n \times n$ alternating sign matrices. Our main 
theorem provides a formula for the number of monotone triangles with prescribed 
bottom row $(k_1,k_2,\ldots,k_n) \in \mathbb{Z}^n$.

\begin{theo}
\label{main}
The number of monotone triangles with $n$ rows and prescribed bottom row $(k_1, k_2, \ldots, k_n)$ 
is given by 
$$
\left( \prod_{1 \le p < q \le n} \left( \id + E_{k_p} \Delta_{k_q} \right) \right) 
\prod_{1 \le i < j \le n} \frac{k_j - k_i}{j-i}, 
$$
where $E_x$ denotes the shift operator, defined by 
$E_x \, p(x) = p(x+1)$, and $\Delta_x:= E_x - \id$ denotes the difference operator.
\end{theo}

In order to understand this formula, there are a few things to remark. 
The product of operators is understood as the composition. 
Moreover note that the shift operators commute, and consequently, it 
does not matter in which order the operators in the product 
$\prod\limits_{1 \le p < q \le n} \left( \id + E_{k_p} \Delta_{k_q} \right)$ are applied.
In order to use this formula to compute the number of monotone triangles with bottom row 
$(k_1,\ldots,k_n)$, one first has to apply the operator
$\prod\limits_{1 \le p < q \le n} \left( \id + E_{x_p} \Delta_{x_q} \right)$ to the polynomial
$\prod\limits_{1 \le i < j \le n} \frac{x_j - x_i}{j-i}$ and then set $x_i=k_i$. Thus, it is not 
so clear how to derive 
\eqref{asmformula} from this formula.

\medskip

Next we discuss the significance of the formula. In the last decades, the enumeration 
of plane partitions, alternating sign matrices and related objects subject to a variety of different 
constraints has attracted a lot of interest. This attraction stems from the fact that now and then these
enumerations lead to an appealing product formula or 
hypergeometric series, 
which is, in spite of their simplicity, pretty hard to prove. At the moment the 
search for these simple product formulas seems to be a bit exhausted.  Therefore, a new challenge is 
the search for possibilities to give enumeration formulas for the vast majority of enumeration problems 
for which there exists no closed formula in a traditional sense. The formula in Theorem~\ref{main} 
contributes to this issue.

Also note that the second product in the formula in Theorem~\ref{main}, i.e. 
$
\prod\limits_{1 \le i < j \le n} \frac{k_j - k_i}{j-i},
$
is the number of semistandard tableaux of shape $(k_n - n, k_{n-1} - n, \ldots, k_1 -1)$ and, 
equivalently, the number of columnstrict plane partitions of this shape 
\cite[p. 375, in (7.105) $q \to 1$]{stanley}. In fact, these objects are in bijection with 
monotone triangles with precsribed bottom row $(k_1,k_2,\ldots,k_n)$
that are strictly increasing in southeast direction \cite[Section 5]{fischer}. 
Thus, our formula once more gives an indication of 
the relation between plane partitions and alternating sign matrices manifested by a number of enumeration 
formulas which show up in both fields, a phenomenon which is not yet well (i.e. bijectively) understood. 
At this point it is also worth mentioning that we can easily rewrite the formula in Theorem~\ref{main} 
such that the second product is the number of semistandard tableaux of shape $(k_n,\ldots,k_1)$. 
\begin{multline*}
\alpha(n;k_1,\ldots,k_n) =
\left( \prod_{1 \le p < q \le n} \left( E^{-1}_{k_q} + E^{-1}_{k_q} E_{k_p} \Delta_{k_q} \right)
\prod_{q=1}^{n} E^{q-1}_{k_{q}} \right) \prod_{1 \le i < j \le n} \frac{k_j - k_i}{j-i}= \\
\left( \prod_{1 \le p < q \le n} \left( \id + E^{-1}_{k_q} \Delta_{k_p} \Delta_{k_q} \right) \right) 
\prod_{1 \le i < j \le n} \frac{k_j - k_i + j -i}{j-i}
\end{multline*}

The method for proving our main theorem can roughly be described as follows. In the first step, we 
introduce a recursion, which relates monotone triangles with $n$ rows to monotone triangles 
with $n-1$ rows. This recursion immediately implies that the enumeration formula  
is a polynomial in $k_1, k_2, \ldots, k_n$. In the next step, we compute the degree of the 
polynomial. Finally, we deduce enough properties of the polynomial in order to 
compute it, where the polynomial's degree determines how much information is in fact needed. This method 
is related to the method for proving polynomial enumeration formulas we have introduced 
in \cite{fischer} and extended in \cite{fischer2}. In the final section we mention some 
further projects around Theorem~\ref{main} we plan to consider next.

\section{The recursion}
\label{recursion}

In the following let $\alpha(n;k_1,\ldots,k_n)$, $n \ge 1$,  denote the number of monotone triangles with $(k_1,\ldots,k_n)$ as bottom row. If we delete the last row of such a monotone triangle we obtain 
a monotone triangle with $n-1$ rows and bottom row, say, $(l_1, l_2, \ldots, l_{n-1})$. By the definition of a monotone triangle we have $k_1 \le l_1 \le k_2 \le l_2 \le \ldots \le k_{n-1} \le l_{n-1} \le k_n$ and $l_i \not= l_{i+1}$. Thus 
\begin{equation}
\label{rec}
\alpha(n;k_1,\ldots,k_n) = 
\sum_{(l_1,\ldots,l_{n-1}) \in \mathbb{Z}^{n-1}, \atop  k_1
\le l_1 \le k_2 \le \ldots \le k_{n-1} \le l_{n-1} \le
k_{n}, l_i \not= l_{i+1}} \alpha(n-1;l_1,\ldots,l_{n-1}).
\end{equation}
We introduce the following abbreviation
$$
\sum_{(l_1,\ldots,l_{n-1}) \in \mathbb{Z}^{n-1}, \atop  k_1 \le
l_1 \le k_2 \le \ldots \le k_{n-1} \le l_{n-1} \le k_{n}, l_i
\not= l_{i+1}} =: \sum^{(k_1,\ldots,k_n)}_{(l_1,\ldots,l_{n-1})}
$$
for $n \ge 2$. This summation operator is well-defined for all 
$(k_1,\ldots,k_n) \in \mathbb{Z}^n$ with $k_1 < k_2 < \ldots < k_n$. We 
extend the definition to arbitrary $(k_1,\ldots,k_n) \in \mathbb{Z}^n$ by 
induction with respect to $n$. If $n=2$ then 
$$
\sum^{(k_1,k_2)}_{(l_1)} A(l_1): =\sum_{l_1=k_1}^{k_2} A(l_1), 
$$
where here and in the following we use the extended definition of the summation over an interval, 
namely,
\begin{equation}
\label{sumext}
\sum_{i=a}^{b} f(i) =
\begin{cases}
f(a) + f(a+1) + \dots + f(b) & \text{if $a \le b$} \\
0 & \text{if $b=a-1$} \\
- f(b+1) - f(b+2) - \dots - f(a-1) & \text{if $b+1 \le a-1 $}
\end{cases}.
\end{equation}
This assures that for any polynomial $p(X)$ over an arbitrary integral domain $I$
containing $\mathbb{Q}$ there exists a unique polynomial $q(X)$ over $I$ such
that $\sum\limits_{x=0}^y p(x) = q(y)$ for all
integers $y$. We usually write $\sum\limits_{x=0}^ y p(x)$ for $q(y)$. (We also use the
analog extended definition for the product symbol.) If $n > 2$ then 
\begin{multline*}
\sum^{(k_1,\ldots,k_n)}_{(l_1,\dots,l_{n-1})} A(l_1,\ldots,l_{n-1}): = \\
\sum^{(k_1,\ldots,k_{n-1})}_{(l_1,\ldots,l_{n-2})} \sum_{l_{n-1}=k_{n-1}+1}^{k_n}
A(l_1,\ldots,l_{n-2},l_{n-1}) +
\sum^{(k_1,\ldots,k_{n-1}-1)}_{(l_1,\ldots,l_{n-2})} A(l_1,\ldots,l_{n-2},k_{n-1}).
\end{multline*}
We renew the definition of $\alpha(n;k_1,\ldots,k_n)$ after this extension by setting $\alpha(1;k_1)=1$ and 
$$
\alpha(n;k_1,\ldots,k_n) = \sum_{(l_1,\ldots,l_{n-1})}^{(k_1,\ldots,k_n)} 
\alpha(n-1;l_1,\ldots, l_{n-1}).
$$
This extends the original function $\alpha(n;k_1,\ldots,k_n)$ to arbitrary 
$(k_1,\ldots,k_n) \in \mathbb{Z}^n$. The recursion implies that $\alpha(n;k_1,\ldots,k_n)$
is a polynomial in $k_1,\ldots,k_n$. We have used this recursion (and a computer) to compute 
$\alpha(n;k_1,\ldots,k_n)$ for $n=1,2,3,4$ and obtain the following
\begin{multline*}
1, 1 - k_1 + k_2, 
\frac{1}{2} ( -3  k_1 + k_1^2 + 2 k_1  k_2 - k_1^2 k_2 - 
        2 k_2^2 + k_1  k_2^2 + 3 k_3 - 4 k_1  k_3 + 
        k_1^2 k_3 + \\  2 k_2 k_3 - k_2^2 k_3 + k_3^2 - 
        k_1 k_3^2 + k_2 k_3^2),
        \frac{1}{12} (20 k_2 + 11 k_1 k_2 - 16 k_1^2 k_2 + 
        3 k_1^3 k_2 + 4 k_1 k_2^2 + 3 k_1^2 k_2^2 - 
        k_1^3 k_2^2 + \\ 4 k_2^3 - 5 k_1 k_2^3 + k_1^2 k_2^3 - 
        20 k_3 + 16 k_1 k_3 - 4 k_1^2 k_3 - 27 k_2 k_3 + 
        9 k_1^2 k_2 k_3 - 2 k_1^3 k_2 k_3 - 
        3 k_1^2 k_2^2 k_3 + k_1^3 k_2^2 k_3 - \\ 
        3 k_2^3 k_3 + 4 k_1 k_2^3 k_3 - k_1^2 k_2^3 k_3 + 
        16 k_1 k_3^2 - 12 k_1^2 k_3^2 + 2 k_1^3 k_3^2 - 
        9 k_1 k_2 k_3^2 + 6 k_1^2 k_2 k_3^2 - 
        k_1^3 k_2 k_3^2 + 9 k_2^2 k_3^2 - \\
        3 k_1 k_2^2 k_3^2 - 3 k_2^3 k_3^2 + 
        k_1 k_2^3 k_3^2 - 4 k_3^3 + 8 k_1 k_3^3 - 
        2 k_1^2 k_3^3 - 3 k_2 k_3^3 - 2 k_1 k_2 k_3^3 + 
        k_1^2 k_2 k_3^3 + 3 k_2^2 k_3^3 - 
        k_1 k_2^2 k_3^3 - \\ 27 k_1 k_4 + 20 k_1^2 k_4 - 
        3 k_1^3 k_4 + 16 k_2 k_4 + 24 k_1 k_2 k_4 - 
        24 k_1^2 k_2 k_4 + 4 k_1^3 k_2 k_4 - 
        16 k_2^2 k_4 + 9 k_1 k_2^2 k_4 + \\
        3 k_1^2 k_2^2 k_4 - k_1^3 k_2^2 k_4 + 
        8 k_2^3 k_4 - 6 k_1 k_2^3 k_4 + k_1^2 k_2^3 k_4 + 
        11 k_3 k_4 - 24 k_1 k_3 k_4 + 15 k_1^2 k_3 k_4 - 
        2 k_1^3 k_3 k_4 - \\ 9 k_2^2 k_3 k_4 + 
        2 k_2^3 k_3 k_4 - 4 k_3^2 k_4 + 9 k_1 k_3^2 k_4 - 
        6 k_1^2 k_3^2 k_4 + k_1^3 k_3^2 k_4 + 
        3 k_2^2 k_3^2 k_4 - k_2^3 k_3^2 k_4 - 
        5 k_3^3 k_4 + \\6 k_1 k_3^3 k_4 - k_1^2 k_3^3 k_4 - 
        4 k_2 k_3^3 k_4 + k_2^2 k_3^3 k_4 - 20 k_1 k_4^2 + 
        9 k_1^2 k_4^2 - k_1^3 k_4^2 + 4 k_2 k_4^2 + 
        15 k_1 k_2 k_4^2 - 9 k_1^2 k_2 k_4^2 + \\
        k_1^3 k_2 k_4^2 - 12 k_2^2 k_4^2 + 
        6 k_1 k_2^2 k_4^2 + 2 k_2^3 k_4^2 - 
        k_1 k_2^3 k_4^2 + 16 k_3 k_4^2 - 
        24 k_1 k_3 k_4^2 + 9 k_1^2 k_3 k_4^2 - 
        k_1^3 k_3 k_4^2 + \\ 9 k_2 k_3 k_4^2 - 
        6 k_2^2 k_3 k_4^2 + k_2^3 k_3 k_4^2 + 
        3 k_3^2 k_4^2 - 3 k_1 k_3^2 k_4^2 + 
        3 k_2 k_3^2 k_4^2 - k_3^3 k_4^2 + 
        k_1 k_3^3 k_4^2 - k_2 k_3^3 k_4^2 - 3 k_1 k_4^3 + \\
        k_1^2 k_4^3 + 2 k_1 k_2 k_4^3 - k_1^2 k_2 k_4^3 - 
        2 k_2^2 k_4^3 + k_1 k_2^2 k_4^3 + 3 k_3 k_4^3 - 
        4 k_1 k_3 k_4^3 + k_1^2 k_3 k_4^3 + 
        2 k_2 k_3 k_4^3 - k_2^2 k_3 k_4^3 + \\ k_3^2 k_4^3 - 
        k_1 k_3^2 k_4^3 + k_2 k_3^2 k_4^3).
\end{multline*}
From this data it is obviously hard to guess a general formula for 
$\alpha(n;k_1,\ldots,k_n)$. However, it seems plausible that the 
degree of $\alpha(n;k_1,\ldots,k_n)$ in $k_i$ is $n-1$. In the following two sections we 
prove that this is indeed true. Note that at first glance the linear growth of the degree 
is quite surprising: suppose $A(l_1,\ldots,l_{n-1})$ is a polynomial of degree no greater 
than $R$ in each of $l_{i-1}$ and $l_{i}$. Then 
\begin{multline*}
\deg_{k_{i}} \left( \sum_{(l_1,\ldots,l_{n-1})}^{(k_1,\ldots,k_n)} A(l_1,\ldots,l_{n-1}) \right)
= \\
\deg_{k_{i}} \left( \sum_{l_{i-1}=k_{i-1}}^{k_i} \sum_{l_{i}=k_i}^{k_{i+1}} A(l_1, \ldots,l_{n-1}) 
- A(l_1,\ldots,l_{i-2}, k_i, k_i, l_{i+1}, \ldots, l_{n-1}) \right) \le 2 R + 2
\end{multline*}
and there exist polynomials $A(l_1,\ldots, l_{n-1})$ such that the upper bound $2R + 2$ is attained. 
Consequently, $\alpha(n;k_1,\ldots, k_n)$ must be of a very specific shape.

\section{Operators related to the recursion}

In this section we define some operators that are fundamental for the study of the recursion 
defined in the previous section. The theory is developed in a bit more generality than it 
is actually needed in order to investigate $\alpha(n;k_1,\ldots,k_n)$.
Recall that the {\it shift operator}, denoted 
by $E_x$, is defined as  $ E_x p(x) = p(x+1).$ Clearly $E_x$ is invertible in the algebra of operators 
of $\mathbb{C}[x]$ and we denote its inverse  by $E^{-1}_x$.
Observe that the shift operators with respect to different variables commute, i.e. $E_{x} E_{y} = E_{y} E_{x}.$ The difference operator $\Delta_x$ is defined as $\Delta_x = E_x - \id$. However, the difference operator
$\Delta_x$ is not invertible since 
it decreases the degree of a polynomial.
If we apply the shift operator or the delta operator to the $i$-th variable of a function, we 
sometimes write $E_i$ or $\Delta_i$, respectively, i.e. 
$\Delta_{k_i} f(k_1,\ldots, k_n) = \Delta_i f(k_1,\ldots, k_n)$. 
Moreover, 
$\Delta_2 f(k_3,k_3,k_3)$, for instance, is shorthand for 
$$
\left. ( \Delta_{l_2} f(l_1,l_2,l_3) ) \right|_{l_1=k_3,l_2=k_3,l_3=k_3}.
$$ 
The {\it swapping operator} $S_{x,y}$ is applicable to functions in (at least) two variables and defined 
as $S_{x,y} f(x,y) = f(y,x).$ 
If we apply $S_{x,y}$ to the $i$-th and $j$-th variable of a function we sometimes write 
$S_{i,j}$.  

In the following we consider rational functions in shift operators. In order to 
guarantee that the inverse  of the denominator always exists, we need the following lemma.

\begin{lem} 
\label{inverse}
Let $p(x_1,\ldots,x_n)$ be a polynomial in $x_1, x_1^{-1}, x_2, x_2^{-1}, \ldots, x_n, x_n^{-1}$ over 
$\mathbb{C}$, and fix an integer $i$, $1 \le i \le n$. Consider the operator
$$
\id + \Delta_{k_i} p(E_{k_1},E_{k_2},\ldots,E_{k_n}) =: \Op
$$
on $\mathbb{C}[k_1,\ldots,k_n]$. Then $\Op$ is invertible and the inverse is 
$$
\Op^{-1} = \sum_{l=0}^{\infty} (-1)^l \Delta_{k_i}^l  p(E_{k_1},E_{k_2},\ldots,E_{k_n})^l, 
$$
where $\Delta_{k_i}^0  p(E_{k_1},E_{k_2},\ldots,E_{k_n})^0 = \id$.
Moreover 
$$
\deg_{k_i} G(k_1,\ldots,k_{n}) = \deg_{k_i} \Op G(k_1,\ldots,k_n) = \deg_{k_i} \Op^{-1} G(k_1,\ldots,k_n).
$$
\end{lem}

{\it Proof.} Let $G(k_1,\ldots,k_n) \in \mathbb{C}[k_1,\ldots,k_n]$. First observe that 
\begin{equation}
\label{opdegree}
\deg_{k_i} G(k_1,\ldots,k_{n}) = \deg_{k_i} \Op G(k_1,\ldots,k_n),
\end{equation}
since $\Delta_{k_i}$ decreases the degree in $k_i$ and $p(E_{k_1},E_{k_2},\ldots,E_{k_n})$
does not increase the degree. It is easy to see that 
$$
F(k_1,\ldots,k_n) = \sum_{l=0}^{\infty} (-1)^l \Delta_{k_i}^l  p(E_{k_1},E_{k_2},\ldots,E_{k_n})^l G(k_1,\ldots,k_n)
$$
is a polynomial with the property that $\Op F = G$. (Observe that the sum is finite since $\Delta^l_{k_i} G(k_1,\ldots,k_n) = 0$ 
if $l > \deg_{k_i} G$.) Assume there is another polynomial $F' \in \mathbb{C}[k_1,\ldots,k_n]$ with the property that $\Op F' = G$. Then 
$\Op H = 0$ with $H=F-F'$. Thus, by \eqref{opdegree}, $\deg_{k_i} H = \deg_{k_i} 0 = - \infty$, a 
contradiction. We obtain $\deg_{k_i} \Op^{-1} G 
= \deg_{k_i} G$ if we apply \eqref{opdegree} to $\Op^{-1} G$. \qed

\medskip

Next we define two operators applicable to polynomials $G(k_1,\ldots,k_n) \in \mathbb{C}[k_1,\ldots,k_n]$. We set 
$$
V_{k_i,k_{j}}= \id + E^{-1}_{k_i} \Delta_{k_i} \Delta_{k_{j}} = 
E^{-1}_{k_i} ( \id + E_{k_{j}} \Delta_{k_i})
$$
and 
$$
T_{k_i,k_{i+1}} = (\id + E_{k_{i+1}} E_{k_i}^{-1} S_{k_i,k_{i+1}} )
\frac{V_{k_i,k_{i+1}}}{V_{k_i,k_{i+1}} + V_{k_{i+1},k_i}}.
$$
By Lemma~\ref{inverse}, the inverse $(V_{k_i,k_{i+1}} + V_{k_{i+1},k_i})^{-1}$  is well-defined. 
The following lemma explains the significance of $T_{k_i,k_{i+1}}$ for the recursion~\eqref{rec}.

\begin{lem}
\label{grad} Let $A(l_1,l_2)$ be a polynomial in $l_1$ and $l_2$
which is of degree at most $R$ in each of $l_1$ and $l_2$. Moreover
assume that $T_{l_1,l_2} A(l_1,l_2)$ is of degree at most $R$ as a
polynomial in $l_1$ and $l_2$, i.e. a linear combination of
monomials $l_1^m l_2^n$ with $m+n \le R$. Then
$$
\sum_{(l_1,l_2)}^{(k_1,k_2,k_3)} A(l_1,l_2) = \sum\limits_{l_1={k_1}}^{k_2} \sum\limits_{l_2={k_2}}^{k_3} A(l_1,l_2) - A(k_2,k_2)
$$
is of degree at most $R+2$ in $k_2$. Moreover, if $T_{l_1,l_2} A(l_1,l_2)=0$
then $\sum\limits_{(l_1,l_2)}^{(k_1,k_2,k_3)} A(l_1,l_2)$ is of degree at most $R+1$ in $k_2$.
\end{lem}

{\it Proof.} We decompose $A(l_1,l_2)= T_{l_1,l_2} A(l_1,l_2) + (\id - T_{l_1,l_2} )
A(l_1, l_2)$. If we define $A^{*}(l_1,l_2)= (\id -  T_{l_1,l_2} )((l_1)_p
(l_2)_q / (p! q!))$, it suffices to show that the degree of
\begin{equation}
\label{expr} \sum\limits_{l_1={k_1}}^{k_2} \sum\limits_{l_2={k_2}}^{k_3}
A^{*}(l_1,l_2)-A^{*}(k_2,k_2)
\end{equation}
in $k_2$ is no greater than $\max(p,q)+1$, where $(x)_p = \prod\limits_{i=0}^{p-1} (x+i)$. Observe that
$$
\id -  T_{l_1,l_2} = \frac{V_{l_2,l_1}}{V_{l_1,l_2} + V_{l_2,l_1}} 
(\id - E_{l_2} E^{-1}_{l_1} S_{l_1,l_2}),
$$
where we use the fact that $S_{l_1,l_2} R(E_{l_1},E_{l_2}) = R(E_{l_2},E_{l_1}) S_{l_1,l_2}$ if 
$R(x,y)$ is a rational function in $x, x^{-1},  y, y^{-1}$.
Lemma~\ref{inverse} implies that
\begin{multline*}
A^{*}(l_1,l_2)= \frac{1}{2} \sum_{i=0}^{\infty} \sum_{j=0}^i \left( -\frac{1}{2} \right)^{i}
\binom{i}{j} \\ \times \left( \frac{(l_1+i-j)_{p-i}
(l_2+j)_{q-i}}{(p-i)! (q-i)!}  - \frac{(l_2+j+1)_{p-i}
(l_1+i-j-1)_{q-i}}{(p-i)! (q-i)!} \right. \\ \left.
+\frac{(l_1+i-j+1)_{p-i-1} (l_2+j)_{q-i-1}}{(p-i-1)! (q-i-1)!} -
\frac{(l_2+j+1)_{p-i-1} (l_1+i-j)_{q-i-1}}{(p-i-1)! (q-i-1)!}
\right).
\end{multline*}
Using the summation formula
\begin{equation}
\label{summation}
\sum_{z=a}^b (z+w)_n = \frac{1}{n+1} ((b+w)_{n+1}-(a-1+w)_{n+1})
\end{equation}
we observe that \eqref{expr} is equal to
\begin{multline*}
\frac{1}{2} \sum_{i=0}^{\infty} \sum_{j=0}^i \left(-\frac{1}{2}\right)^{i}
\binom{i}{j} \\ \times \left( - \binom{{k_2}+p-j}{p-i+1}
\binom{{k_2}+q+j-i-1}{q-i+1} + \binom{{k_2}+q-j-1}{q-i+1}
\binom{{k_2}+p+j-i}{p-i+1} \right. \\
 - \binom{{k_2}+p-j}{p-i}
\binom{{k_2}+q+j-i-2}{q-i} + \binom{{k_2}+q-j-1}{q-i}
\binom{{k_2}+p+j-i-1}{p-i} \\
 - \binom{{k_2}+p-j-1}{p-i}
\binom{{k_2}+q+j-i-1}{q-i} + \binom{{k_2}+q-j-2}{q-i} \binom{{k_2}+p+j-i}{p-i}
\\
\left.  - \binom{{k_2}+p-j-1}{p-i-1} \binom{{k_2}+q+j-i-2}{q-i-1} +
\binom{{k_2}+q-j-2}{q-i-1} \binom{{k_2}+p+j-i-1}{p-i-1} \right) \\
 + R(k_1, k_2, k_3),
\end{multline*}
where $R(k_1, k_2, k_3)$ is a polynomial in $k_1,k_2,k_3$ of degree no greater than
$\max(p,q)+1$ in $k_2$. If we replace $j$ by $i-j$ in every other product
of two binomial coefficients we see that this expression
simplifies to $R(k_1, k_2, k_3)$ and the lemma is proved. \qed

\medskip

In order to use Lemma~\ref{grad} to compute the degree of $\sum\limits_{(l_1,l_2)}^{(k_1,k_2,k_3)} A(l_1,l_2)$ 
in $k_2$, one has to compute the degree of $T_{l_1,l_2} A(l_1,l_2)$ in $l_1$ and $l_2$. However, the 
operator $T_{l_1,l_2}$ is complicated, and thus it is convenient to consider a simplified version of 
$T_{l_1,l_2}$ for this purpose, which is obtained by multiplication with an operator that preserves the degree.
\begin{multline*}
T'_{k_i,k_{i+1}}:= E_{k_i} (V_{k_i,k_{i+1}} + V_{k_{i+1}, k_{i}} ) T_{k_i,k_{i+1}} = \\
(\id + S_{k_i,k_{i+1}}) E_{k_i} V_{k_i,k_{i+1}} = 
(\id + S_{k_i,k_{i+1}}) (\id + E_{k_{i+1}} \Delta_{k_i}) 
\end{multline*}
Observe that $\deg_{k_i,k_{i+1}} T_{k_i,k_{i+1}} G(k_1,\ldots,k_n) = \deg_{k_i,k_{i+1}} T'_{k_i,k_{i+1}} G(k_1,\ldots,k_n)$, since 
$$V_{k_i,k_{i+1}}+V_{k_{i+1},k_i}= 2 \id + (E^{-1}_{k_i} + E^{-1}_{k_{i+1}}) \Delta_{k_i} \Delta_{k_{i+1}}$$ 
and $\Delta_{k_i} \Delta_{k_{i+1}}$ decreases the degree of a polynomial in 
$k_{i}$ and $k_{i+1}$. In particular, $T_{k_i,k_{i+1}} G(k_1,\ldots,k_n)=0$ if and only if 
$T'_{k_i,k_{i+1}} G(k_1,\ldots,k_n)=0$. 

\section{The fundamental lemma}

Suppose $A(l_1,\ldots,l_{n})$ is a function on $\mathbb{Z}^{n}$. In this section we prove a lemma that 
expresses 
$$
T'_{k_i,k_{i+1}} \left( \sum\limits_{(l_1,\ldots,l_{n})}^{(k_1,\ldots,k_{n+1})} A(l_1,\ldots,l_n) \right)
(k_1,\ldots,k_{n+1})
$$
in terms of $T'_{l_{i-1},l_{i}} A(l_1,\ldots,l_n)$ and $T'_{l_i,l_{i+1}} A(l_1,\dots,l_n)$. In particular,
this shows that if $T'_{l_i,l_{i+1}} A(l_1,\dots,l_n)=0$ for all $i=1,\ldots, n-1$ then 
$$T'_{k_i,k_{i+1}} \left(  \sum\limits_{(l_1,\ldots,l_{n})}^{(k_1,\ldots,k_{n+1})} A(l_1,\ldots,l_n) \right)
(k_1,\ldots,k_{n+1}) =0$$ for all $i=1,\ldots, n$.

\begin{lem} 
\label{fundASM}
Let $f(k_1,k_2,k_3)$ be a function from $\mathbb{Z}^3$ to 
$\mathbb{C}$ and define
$$g(k_1,k_2,k_3,k_4): = \sum_{(l_1,l_2,l_3)}^{(k_1,k_2,k_3,k_4)} f(l_1,l_2,l_3). 
$$
Then 
\begin{multline*}
T'_{2,3} \, g(k_1,k_2,k_3,k_4) = \\
- \frac{1}{2} \left(  \sum_{l_1=k_2+1}^{k_3} \sum_{l_2=k_2+1}^{k_3} \sum_{l_3=k_2}^{k_4} T'_{1,2} \, f(l_1,l_2,l_3) + \sum_{l_1=k_1}^{k_2+1} \sum_{l_2=k_2}^{k_3-1} \sum_{l_3=k_2}^{k_3-1} 
T'_{2,3} \, f(l_1,l_2,l_3) \right) \\
+ \frac{1}{2} \left( \sum_{l_1=k_2}^{k_3-1} \sum_{l_2=k_2}^{k_3-1} \Delta_2  (\id + E_1) 
T'_{1,2} \, f(l_1,l_2,k_2) 
- \sum_{l_2=k_2}^{k_3-1} \sum_{l_3=k_2}^{k_3-1} \Delta_2 (\id + E_3) 
T'_{2,3} \, f(k_2+1,l_2,l_3) \right) \\
+ \frac{1}{2} \Big( T'_{1,2} \, f(k_2,k_2,k_2+1) - T'_{1,2} \, f(k_2,k_2,k_3+1) + T'_{2,3} \, f(k_2,k_2,k_2)
- T'_{2,3} \, f(k_3,k_2,k_2) \Big) \\
- T'_{1,2} \, f(k_2,k_3,k_2+1) - T'_{2,3} \, f(k_2,k_2,k_3). 
\end{multline*}
Moreover, for a function $h(l_1, l_2)$ on $\mathbb{Z}^2$,
$$
T'_{1,2} \left( \sum_{(l_1,l_2)}^{(k_1,k_2,k_3)} h(l_1,l_2) \right)(k_1,k_2,k_3) = - \frac{1}{2} \sum_{l_1=k_1}^{k_2-1} \sum_{l_2=k_1}^{k_2-1} T'_{1,2} \, h(l_1,l_2).
$$
\end{lem}  

{\it Proof.} 
We only sketch the proof of the first formula since the proof of the second formula is easy. Observe that in this 
formula $T'_{1,2} f(k_2,k_2,k_{2}+1)$, for instance, is shorthand for 
$$
\left. T'_{l_1,l_2} f(l_1,l_2,k_{2}+1) \right|_{l_1=k_2,l_2=k_2}.
$$
By definition
$$
g(k_1,k_2,k_3,k_4) = 
\sum_{l_1=k_1}^{k_2} \sum_{l_2=k_2}^{k_3} \sum_{l_3=k_3}^{k_4} f(l_1,l_2,l_3) - 
\sum_{l_3=k_3}^{k_4} f(k_2,k_2,l_3) - \sum_{l_1=k_1}^{k_2} f(l_1,k_3,k_3).
$$
It is easy to see that 
$$
\Delta_2 g(k_1,k_2,k_3,k_4) = 
- \sum_{l_1=k_1}^{k_2-1} \sum_{l_3=k_3}^{k_4} f(l_1, k_2, l_3) + 
\sum_{l_2=k_2+2}^{k_3} \sum_{l_3=k_3}^{k_4} f(k_2+1,l_2,l_3) - f(k_2+1,k_3,k_3).
$$
This implies that 
\begin{multline*}
(\id + \Delta_2 E_{3}) g(k_1,k_2,k_3,k_4) = 
\sum_{l_1=k_1}^{k_2} \sum_{l_2=k_2}^{k_3} \sum_{l_3=k_3}^{k_4} f(l_1,l_2,l_3) \\  
- \sum_{l_3=k_3}^{k_4} f(k_2,k_2,l_3) - \sum_{l_1=k_1}^{k_2} f(l_1,k_3,k_3) - f(k_2+1,k_3+1,k_3+1) \\ -
\sum_{l_1=k_1}^{k_2-1} \sum_{l_3=k_3+1}^{k_4} f(l_1,k_2,l_3) 
+ \sum_{l_2=k_2+2}^{k_3+1} \sum_{l_3=k_3+1}^{k_4} f(k_2+1,l_2,l_3).
\end{multline*}
Next we want to apply the operator $(\id+S_{2,3})$ .
Observe that 
\begin{multline*}
(\id+S_{2,3}) \left( \sum_{l_1=k_1}^{k_2} \sum_{l_2=k_2}^{k_3} \sum_{l_3=k_3}^{k_4} f(l_1,l_2,l_3) 
\right)(k_1,k_2,k_3,k_4) \\ = 
\sum_{l_1=k_1}^{k_2} \sum_{l_2=k_2}^{k_3} \sum_{l_3=k_3}^{k_4} f(l_1,l_2,l_3) -
\sum_{l_1=k_1}^{k_3} \sum_{l_2=k_2}^{k_3} \sum_{l_3=k_2}^{k_4} f(l_1,l_2,l_3)  \\
= - \sum_{l_1=k_2+1}^{k_3} \sum_{l_2=k_2+1}^{k_3} \sum_{l_3=k_2}^{k_4} f(l_1,l_2,l_3)
- \sum_{l_1=k_1}^{k_2} \sum_{l_2=k_2}^{k_3-1} \sum_{l_3=k_2}^{k_3-1} f(l_1,l_2,l_3) 
+ \sum_{l_1=k_1}^{k_2} \sum_{l_3=k_2}^{k_4} f(l_1,k_2,l_3) \\
+ \sum_{l_1=k_1}^{k_3} \sum_{l_3=k_3}^{k_4} f(l_1,k_3,l_3) 
+ \sum_{l_1=k_2+1}^{k_3} \sum_{l_3=k_2}^{k_3-1} f(l_1,k_3,l_3)  \\
= - \frac{1}{2} \left( \sum_{l_1=k_2+1}^{k_3} \sum_{l_2=k_2+1}^{k_3} \sum_{l_3=k_2}^{k_4} (\id+S_{1,2}) f(l_1,l_2,l_3) + \sum_{l_1=k_1}^{k_2} \sum_{l_2=k_2}^{k_3-1} \sum_{l_3=k_2}^{k_3-1} 
(\id+S_{2,3}) f(l_1,l_2,l_3)  \right) \\
+ \sum_{l_1=k_1}^{k_2} \sum_{l_3=k_2}^{k_4} f(l_1,k_2,l_3) 
+ \sum_{l_1=k_1}^{k_3} \sum_{l_3=k_3}^{k_4} f(l_1,k_3,l_3) 
+ \sum_{l_1=k_2+1}^{k_3} \sum_{l_3=k_2}^{k_3-1} f(l_1,k_3,l_3).
\end{multline*}
Therefore, we have 
\begin{multline*}
(\id+S_{2,3})(\id + \Delta_2 E_3) g(k_1,k_2,k_3,k_4)  \\
= -\frac{1}{2} \left( 
\sum_{l_1=k_2+1}^{k_3} \sum_{l_2=k_2+1}^{k_3} \sum_{l_3=k_2}^{k_4} 
T'_{1,2} \, f(l_1,l_2,l_3)
 +
\sum_{l_1=k_1}^{k_2+1} \sum_{l_2=k_2}^{k_3-1} \sum_{l_3=k_2}^{k_3-1} 
T'_{2,3} \, f(l_1,l_2,l_3)   \right) \\
+ \sum_{l_3=k_2}^{k_3} f(k_2+1,k_2+1,l_3) 
- \sum_{l_3=k_2}^{k_3+1} f(k_2+1,k_3+1,l_3) 
- \sum_{l_2=k_2}^{k_3-1} f(k_2+1,l_2,k_3) \\ 
+ \sum_{l_2=k_2+2}^{k_3+1} f(k_3+1,l_2,k_2)  
+ \sum_{l_1=k_2+2}^{k_3-1} f(l_1,k_3,k_2) - 
\sum_{l_1=k_2+2}^{k_3} f(l_1,k_2,k_2) \\
- f(k_2,k_2,k_3) - f(k_3+1,k_2+1,k_2+1).
\end{multline*}
Finally check that the right-hand-side of this equation is equal to the right-hand-side in 
the statement of the lemma.
\qed

\medskip

This proves the statement preceding the lemma for $n=2,3$. It can
easily be extended to general $n$ by deriving a merging rule for the recursion \eqref{rec}. 
For this purpose we need another operator. Let $f(x,z)$ be a 
function on $\mathbb{Z}^2$. Then the operator $I_{x,z}^y$ transforms 
$f(x,z)$ into a function on $\mathbb{Z}$ by 
$$
I_{x,z}^{y} f(x,z) := f(y-1,y) + f(y,y+1) - f(y-1,y+1) = \left. V_{x,z} f(x,z) \right|_{x=y, z=y}.
$$
With this definition we have
\begin{equation}
\label{merge}
\sum_{(l_1,\ldots,l_{n-1})}^{(k_1,\ldots,k_n)} A(l_1,\ldots,l_{n-1}) =
I_{k'_i,k''_i}^{k_i} 
\sum_{(l_1,\ldots,l_{i-1})}^{(k_1,\ldots,k_{i-1},k'_i)} 
\sum_{(l_i,\ldots,l_{n-1})}^{(k''_i,k_{i+1},\ldots, k_n)} 
A(l_1,\ldots, l_n).
\end{equation}
Fix a function $A(l_1,\ldots,l_n)$ on $\mathbb{Z}^n$ and an $i$ with $2 \le i \le n-1$.
Let 
$$
A_{x,y}'(l_1, \ldots, l_{i-2}, k_i, k_{i+1}, l_{i+2}, \ldots, l_n) =
\sum_{(l_{i-1},l_i,l_{i+1})}^{(x,k_i,k_{i+1},y)} 
A(l_1,\ldots, l_{n})
$$ 
and 
\begin{multline*}
A''_{w,x,y,z}(k_1,\ldots,k_{i-2}, k_i, k_{i+1}, k_{i+3}, \ldots, k_{n+1}) = 
\sum_{(l_1,\ldots,l_{i-2})}^{(k_1,\ldots,k_{i-2},w)} 
\sum_{(l_{i+2},\ldots,l_n)}^{(z,k_{i+3},\ldots,k_{n+1})}  A_{x,y}'(l_1,\ldots,l_n).
\end{multline*}
Then, by \eqref{merge},  
$$
\sum_{(l_1,\ldots,l_n)}^{(k_1,\ldots,k_{n+1})} A(l_1,\ldots,l_n) = 
I_{w,x}^{k_{i-1}} I_{y,z}^{k_{i+2}} A''_{w,x,y,z}(k_1,\ldots,k_{i-2},k_i,k_{i+1},k_{i+3}, \ldots, k_{n+1}).
$$
Define \begin{multline*}
A_{x,y}^{*}(l_1, \ldots, l_{i-2}, k_i, k_{i+1}, l_{i+2}, \ldots, l_n) = \\
- \frac{1}{2} \left( \sum_{l_{i-1}=k_i+1}^{k_{i+1}} \sum_{l_i=k_i+1}^{k_{i+1}} 
\sum_{l_{i+1}=k_i}^{y} T'_{i-1,i} A(l_1,\ldots, l_n) + \sum_{l_{i-1}=x}^{k_i+1} \sum_{l_i=k_i}^{k_{i+1}-1} 
\sum_{l_{i+1}=k_i}^{k_{i+1}-1} T'_{i,i+1} A(l_1,\ldots,l_n) \right) \\
+ \frac{1}{2} \left( \sum_{l_{i-1}=k_i}^{k_{i+1}-1} \sum_{l_i=k_{i}}^{k_{i+1}-1}  \Delta_i  (\id + E_{i-1}) 
T'_{i-1,i} A(l_1,\ldots, l_i, k_{i}, l_{i+2}, \ldots, l_n)  \right. \\ \left. 
- \sum_{l_i=k_i}^{k_{i+1}-1} \sum_{l_{i+1}=k_i}^{k_{i+1}-1}  \Delta_i (\id + E_{i+1}) 
T'_{i,i+1} A(l_1,\ldots, l_{i-2}, k_{i}+1, l_i, \ldots, l_n) \right)  \\
\frac{1}{2} \Big( 
T'_{i-1,i} A(\ldots, l_{i-2}, k_i, k_i, k_i +1, l_{i+2}, \ldots) -
T'_{i-1,i} A(\ldots, l_{i-2}, k_i, k_i, k_{i+1} +1, l_{i+2}, \ldots)  \\
  +
T'_{i,i+1} A(\ldots, l_{i-2}, k_i, k_i, k_{i}, l_{i+2}, \ldots) -
T'_{i,i+1} A(\ldots, l_{i-2}, k_{i+1}, k_i, k_{i}, l_{i+2}, \ldots) \Big) \\
- T'_{i-1,i} A(\ldots, l_{i-2}, k_i, k_{i+1}, k_i +1, l_{i+2}, \ldots)
- T'_{i,i+1} A(\ldots, l_{i-2}, k_i, k_{i}, k_{i+1}, l_{i+2}, \ldots)
\end{multline*}
and 
\begin{multline*}
A^{**}_{w,x,y,z}(k_1,\ldots,k_{i-2}, k_i, k_{i+1}, k_{i+3}, \ldots, k_{n+1}) = 
\sum_{(l_1,\ldots,l_{i-2})}^{(k_1,\ldots,k_{i-2},w)} 
\sum_{(l_{i+2},\ldots,l_n)}^{(z,k_{i+3},\ldots,k_{n+1})} A_{x,y}^{*}(l_1,\ldots,l_n).
\end{multline*}
Then, by the first formula in Lemma~\ref{fundASM}, we have  
\begin{multline}
\label{connect}
T'_{k_i,k_{i+1}} \left( \sum_{(l_1,\ldots,l_n)}^{(k_1,\ldots,k_{n+1})} 
A(l_1,\ldots,l_n) \right) (k_1,\ldots,k_{n+1}) \\ = 
I_{w,x}^{k_{i-1}} I_{y,z}^{k_{i+2}} A^{**}_{w,x,y,z}(k_1,\ldots,k_{i-2},k_i,k_{i+1},k_{i+3},\ldots, k_{n+1}).
\end{multline}
If we use the second formula in Lemma~\ref{fundASM}, we obtain a similar formula for the 
case $i=1$. By symmetry an analog formula follows for $i=n$. These formulas imply the following corollary.

\begin{cor} Suppose $A(l_1,\ldots,l_n)$ is a function on $\mathbb{Z}^n$ with 
$T'_{l_i,l_{i+1}} A(l_1,\ldots,l_n)=0$ for all $i$, $1 \le i < n$. Then 
$$
T'_{k_i,k_{i+1}} \left( \sum_{(l_1,\ldots,l_n)}^{(k_1,\dots,k_{n+1})} A(l_1,\ldots,l_n) \right)
(k_1, \ldots, k_{n+1}) =0
$$
for all $i$, $1 \le i \le n$.
\end{cor}

We come back to $\alpha(n;k_1,\ldots,k_n)$.
By induction with respect to $n$ we conclude that $T'_{k_i,k_{i+1}} \alpha(n;k_1,\dots,k_n) = 0$ for 
all $i$, $1 \le i < n$, if $n \ge 2$. (Note that $\alpha(2;k_1,k_2)=k_2-k_1+1$.) Thus $T_{k_i,k_{i+1}} \alpha(n;k_1,\ldots,k_n)=0$ for all $i$. Therefore, by Lemma~\ref{grad}
and by induction with respect to $n$, the polynomial $\alpha(n;k_1,\ldots,k_n)$ is 
of degree no greater than $n-1$ in every $k_i$. 

\section{Proof of the theorem}

In the previous two sections we have seen that the fact that $T'_{k_i,k_{i+1}} \alpha(n;k_1,\ldots,k_n)=0$ 
for all $i$ is fundamental for the computation of the polynomial's degree. In this section, however, we 
demonstrate that this property already determines $\alpha(n;k_1,\ldots,k_n)$ up to a multiplicative 
constant. Observe that $T'_{k_i,k_{i+1}} A(k_1,\dots,k_n)=0$ is equivalent with the fact that 
$(\id + E_{k_{i+1}} \Delta_{k_i} ) A(k_1,\ldots,k_n)$ is antisymmetric in $k_i$ and 
$k_{i+1}$. In the following lemma we characterize functions $A(k_1,\ldots, k_n)$ with the 
property that $$(\id + E_{k_{i+1}} \Delta_{k_i} ) A(k_1,\ldots,k_n)$$ is antisymmetric in $k_i$ and $k_{i+1}$
for all $i$.

\begin{lem} Let $A(k_1,\ldots,k_n)$ be a polynomial in $(k_1,\ldots,k_n)$. Then 
$$
(\id + E_{k_{i+1}} \Delta_{k_i} ) A(k_1,\ldots,k_n)
$$
is antisymmetric in $k_i$ and $k_{i+1}$ for all $i$, $1 \le i \le n-1$, if and only if 
$$
\left( \prod_{1 \le p < q \le n} (\id + E_{k_{q}} \Delta_{k_p}) \right) A(k_1,\ldots,k_n)
$$
is antisymmetric in $k_1,\ldots,k_n$.
\end{lem}

{\it Proof.} First assume that $(\id  + E_{k_{i+1}} \Delta_{k_i}) A(k_1,\ldots,k_n)$
is antisymmetric in $k_i$ and $k_{i+1}$ for all $i$. We have to show that 
$$
(\id + S_{k_i,k_{i+1}}) \left( \prod_{1 \le p < q \le n} (\id + E_{k_q} \Delta_{k_p} ) \right) A(k_1,\ldots,k_n) =0
$$
for all $i$. For this purpose observe that 
\begin{multline*}
(\id + S_{k_i,k_{i+1}}) \left( \prod_{1 \le p < q \le n, (p,q) \not= (i,i+1)} (\id  + E_{k_q} \Delta_{k_p}) \right)
(\id+ E_{k_{i+1}} \Delta_{k_i} )  A(k_1,\ldots,k_n) \\
= \left( \prod_{1 \le p < q \le n, (p,q) \not= (i,i+1)} (\id+ E_{k_q} \Delta_{k_p} ) \right) (\id + S_{k_i,k_{i+1}}) 
(\id  + E_{k_{i+1}} \Delta_{k_i} )  A(k_1,\ldots,k_n) =0,
\end{multline*} 
because 
\begin{multline*}
\prod\limits_{1 \le p < q \le n, (p,q) \not= (i,i+1)} (\id + E_{k_q} \Delta_{k_p} )
= \left( \prod\limits_{1 \le p < q \le n, \atop p,q \notin \{i,i+1\}} (\id + E_{k_q} \Delta_{k_p} ) \right) 
\left( \prod_{i+1 < q \le n} (\id + E_{k_q} \Delta_{k_i}) \right) \\
\times
\left( \prod_{i+1 < q \le n} (\id + E_{k_q} \Delta_{k_{i+1}}) \right) 
\left( \prod_{1 \le p < i} (\id + E_{k_i} \Delta_{k_p}) \right)
\left( \prod_{1 \le p < i} (\id + E_{k_{i+1}} \Delta_{k_p}) \right)
\end{multline*}
is symmetric in 
$k_i$ and $k_{i+1}$. Conversely, assume that 
$$
\left( \prod_{1 \le p < q \le n} (\id + E_{k_q} \Delta_{k_p} ) \right) A(k_1,\ldots,k_n)
$$
is antisymmetric in $k_1,\ldots,k_n$. Consequently, 
$$
\left( \prod_{1 \le p < q \le n, (p,q) \not= (i,i+1)} (\id + E_{k_q} \Delta_{k_p} ) \right) 
(\id + S_{k_i,k_{i+1}}) 
(\id + E_{k_{i+1}}  \Delta_{k_i} )  A(k_1,\ldots,k_n) =0,
$$
for all $i$, $1 \le i \le n-1$. By Lemma~\ref{inverse} the operator $\prod\limits_{1 \le p < q \le n, \atop (p,q) \not= (i,i+1)} (\id  +  E_{k_q} \Delta_{k_p} )$ is invertible, and therefore
$
(\id + S_{k_i,k_{i+1}}) 
(\id + E_{k_{i+1}} \Delta_{k_i} )  A(k_1,\ldots,k_n) =0. \qed
$

\medskip

Using this lemma we see that 
\begin{equation}
\label{stern}
\left( \prod_{1 \le p < q \le n} (\id + E_{k_q} \Delta_{k_p} ) \right) \alpha(k_1,\ldots,k_n)
\end{equation}
is an antisymmetric polynomial in $k_1,\ldots,k_n$. A product of shift operators does 
not increase a polynomial's degree, and thus the degree of \eqref{stern} in every $k_i$ is no 
greater than $n-1$.  Every antisymmetric function in $k_1,\ldots, k_n$ is a multiple of 
$\prod\limits_{1 \le i < j \le n} (k_j - k_i)$,
and since this product is of degree $n-1$ in every $k_i$, the expression in \eqref{stern}
is equal to $ C  \prod\limits_{1 \le i < j \le n} (k_j - k_i)$, where $C$ is a rational constant. 
By Lemma~\ref{inverse} $\prod\limits_{1 \le p < q \le n} (\id + E_{k_q} \Delta_{k_p} )$ is invertible, and
therefore
$$
\alpha(n;k_1,\ldots,k_n) =  \left( \prod_{1 \le p < q \le n} \frac{1}{\id  + E_{k_q} \Delta_{k_p} }
\right) C 
\prod_{1 \le i < j \le n} (k_j - k_i).
$$
We compute the constant $C$. We expand $\alpha(n;k_1,\ldots,k_n)$ with respect to the 
basis $\prod\limits_{i=1}^{n} (k_i)_{m_i}$ and consider the (non-zero) coefficient of the basis element 
with maximal $(m_n,m_{n-1},\ldots,m_1)$ in lexicographic order. We show by induction with respect to $n$ that 
$(m_n,m_{n-1},\ldots,m_1)=(n-1,n-2,\ldots,1,0)$ and that the coefficient is 
$\prod\limits_{i=1}^{n} \frac{1}{(i-1)!}$. Assume that the assertion is true for $n-1$.
A careful analysis of the definition of $\sum\limits_{(l_1,\ldots,l_{n-1})}^{(k_1,\ldots,k_n)}$ shows that 
the ``maximal'' basis element of $\alpha(n;k_1,\ldots,k_n)$ with respect to the  
lexicographic order is the ``maximal'' basis element of 
$$
\sum_{l_1=k_1}^{k_2} \sum_{l_2=k_2+1}^{k_3} \ldots \sum_{l_{n-1}=k_{n-1}+1}^{k_n} 
\prod_{i=1}^{n-1} \frac{(l_i)_{i-1}}{(i-1)!}.
$$
The assertion follows and thus 
$
C= \prod\limits_{i=1}^{n} \frac{1}{(i-1)!} = \prod\limits_{1 \le i < j \le n} \frac{1}{j-i}.
$
We obtain the following Theorem.

\begin{theo}
\label{maininverse}
The number of monotone triangles with $n$ rows and prescribed bottom row $(k_1, k_2, \ldots, k_n)$ 
is equal to
$$
\left( \prod_{1 \le p < q \le n} \frac{1}{\id + E_{k_q} \Delta_{k_p}} \right) 
\prod_{1 \le i < j \le n} \frac{k_j - k_i}{j-i}.
$$
\end{theo}

By the formula for the geometric series, the inverse of the operator
$\id + E_{k_q} \Delta_{k_p}$ appearing in this formula is equal to 
$$
\sum_{l=0}^{\infty} (-1)^l E^{l}_{k_q} \Delta^{l}_{k_p}.
$$
This follows from the proof of Lemma~\ref{inverse}. However, it is 
also possible to give a similar formula for $\alpha(n;k_1,\ldots,k_n)$ 
which does not involve inverses of operators. In order to derive it, we need 
the following lemma.

\begin{lem} 
\label{sym}
Let $P(X_1,\ldots,X_n)$ be a polynomial in $(X_1,\ldots,X_n)$ over $\mathbb{C}$
which is symmetric in $(X_1,\ldots, X_n)$. Then 
$$
P(E_{k_1},\ldots, E_{k_n}) \prod_{1 \le i < j \le n} \frac{k_j - k_i}{j-i} = 
P(1,\ldots,1) \cdot \prod_{1 \le i < j \le n} \frac{k_j - k_i}{j-i}.
$$
\end{lem}

{\it Proof.} Let $(m_1,\ldots,m_n) \in \mathbb{Z}^n$ be with $m_i \ge 0$ for all $i$ and $m_i \not= 0$ for at 
least one $i$. It suffices to show that 
$$
\sum_{\pi \in {\mathcal S}_n} \Delta^{m_{\pi(1)}}_{k_1} \Delta^{m_{\pi(2)}}_{k_2} \dots 
\Delta^{m_{\pi(n)}}_{k_n}  \prod_{1 \le i < j \le n} \frac{k_j - k_i}{j-i} = 0. $$
By the Vandermonde determinant evaluation, we have 
$$
\prod_{1 \le i < j \le n} \frac{k_j - k_i}{j-i} = \det_{1 \le i, j \le n} \left( \binom{k_i}{j-1} \right).
$$
Therefore, it suffices to show that 
$$
\sum_{\pi, \sigma \in {\mathcal S}_n} \sgn \sigma \binom{k_1}{\sigma(1)-m_{\pi(1)}-1}
\binom{k_2}{\sigma(2)-m_{\pi(2)}-1} \dots \binom{k_n}{\sigma(n)-m_{\pi(n)}-1}=0.
$$
If, for fixed $\pi, \sigma \in {\mathcal S}_n$,  there exists an $i$ with $\sigma(i)-m_{\pi(i)}-1 < 0$ then the corresponding 
summand vanishes. We define a sign reversing involution on the set of non-zero summands. 
Fix $\pi, \sigma \in {\mathcal S}_n$ such that the summand corresponding to $\pi$ and $\sigma$ 
does not vanish. Consequently, $\{\sigma(1)-m_{\pi(1)}-1,\sigma(2)-m_{\pi(2)}-1,\ldots,\sigma(n)-m_{\pi(n)}-1\} \subseteq \{0,1,\ldots,n-1\}$ and since $(m_1,\ldots,m_n) \not= (0,\ldots,0)$ there are $i, j$, $1 \le i < j \le n$, 
with $\sigma(i)-m_{\pi(i)}-1=\sigma(j)-m_{\pi(j)}-1$. Among all pairs $(i,j)$ with this property, let $(i',j')$ be the pair which is 
minimal with respect to the lexicographic order. Then the summand corresponding to 
$\pi \circ (i',j')$ and $\sigma \circ (i',j')$ is the negativ of the summand corresponding to $\pi$ and $\sigma$. \qed

\medskip

Observe that $\prod\limits_{1 \le p, q \le n} (1 + X_q ( X_p - 1))$ is symmetric in $(X_1,\ldots, X_n)$. Thus, by 
Lemma~\ref{sym}, 
$$
\prod_{1 \le p, q \le n} \left( \id + E_{k_q} \Delta_{k_p} \right) 
\prod_{1 \le i < j \le n} \frac{k_j - k_i}{j-i} =  \prod_{1 \le i < j \le n} \frac{k_j - k_i}{j-i}.
$$
Therefore, by Theorem~\ref{maininverse},  
\begin{multline*}
\alpha(n;k_1,\ldots,k_n) = \left( \prod_{1 \le p < q \le n} \frac{1}{\id + E_{k_q} \Delta_{k_p}} \right) 
\prod_{1 \le i < j \le n} \frac{k_j - k_i}{j-i} \\
= \left( \prod_{1 \le p < q \le n} \frac{1}{\id + E_{k_q} \Delta_{k_p}} \right)
\left( \prod_{1 \le p, q \le n} \left( \id + E_{k_q} \Delta_{k_p} \right) \right)
\prod_{1 \le i < j \le n} \frac{k_j - k_i}{j-i} \\
= \left( \prod_{1 \le p < q \le n} \left( \id + E_{k_p} \Delta_{k_q} \right) \right)
\prod_{1 \le i < j \le n} \frac{k_j - k_i}{j-i}
\end{multline*}
and this completes the proof of Theorem~\ref{main}.

\section{Some further projects}

In this section we list some further projects around 
the formula given in Theorem~\ref{main} we plan to pursue.  

\begin{enumerate}

\item A natural question to ask is whether it is possible to derive the formula 
for the number of $n \times n$ alternating sign matrices \eqref{asmformula} from Theorem~\ref{main}, 
i.e. to show that
$$
\left. \left[ \left( \prod_{1 \le p < q \le n} \left( \id + E_{k_p} \Delta_{k_q} \right) \right) 
\prod_{1 \le i < j \le n} \frac{k_j - k_i}{j-i} \right] \right|_{(k_1,k_2,\ldots,k_n)=(1,2,\ldots,n)} 
= \prod_{j=1}^{n} \frac{(3j-2)!}{(n+j-1)!} 
$$
More generally, one could try to reprove the refined alternating sign matrix 
theorem \cite{zeilberger2}, which states that the number of $n \times n$ 
alternating sign matrices in which the unique $1$ in the top row is in the 
$k$-th column is given by 
\begin{equation}
\label{refined}
\frac{(k)_{n-1}
(1+n-k)_{n-1}}{(n-1)!} \prod\limits_{j=1}^{n-1} \frac{(3j-2)!}{(n+j-1)!}.
\end{equation}
An analysis of the correspondence between alternating sign matrices and 
monotone triangles shows that $\alpha(n-1;1,2,\ldots,k-1,k+1,\ldots,n)$ is 
the number of $n \times n$ alternating sign matrices in which the unique 
$1$ in the bottom row is in the $k$-th column and this is by symmetry 
equal to \eqref{refined}. This could be a 
consequence of an even more general theorem: computer experiments 
suggest that there are other $(k_1,k_2,\ldots,k_n) \in \mathbb{Z}^n$  
``near'' $(1,2,\ldots,n)$ for which $\alpha(n;k_1,\ldots,k_{n})$ has only small 
prime factors. Small prime factors are an indication for a 
simple product formula. A similar phenomenon can be observed for some  
$(k_1,k_2,\ldots,k_n) \in \mathbb{Z}^n$ ``near'' $(1,3,\ldots,2n-1)$.
It is not too hard to see that $\alpha(n;1,3,\ldots,2n-1)$ is 
the number of $(2n+1) \times (2n+1)$ alternating sign matrices, which are symmetric 
with respect to the reflection along the vertical axis. Kuperberg \cite{kuperberg2} showed that the 
number of these objects is 
$$
\frac{n!}{(2n)! 2^n} \prod_{j=1}^n \frac{(6j-2)!}{(2n+2j-1)!} .
$$

\item Let $\beta(n;k_1,\ldots,k_n)$ denote the number of monotone triangles with 
prescribed bottom row $(k_1,\ldots, k_n)$ that are strictly increasing in southeast 
direction. With this notation, Theorem~\ref{main} states that
\begin{equation}
\label{equi}
\alpha(n;k_1,\ldots,k_n) = \left( \prod_{1 \le p < q \le n} (\id + E_{k_p} \Delta_{k_q}) \right) 
\beta(n;k_1,\ldots,k_n).
\end{equation}
It would be interesting to find a bijective proof of this formula 
in the following sense: if we expand the product of operators on the 
left hand side we obtain a sum of expressions of the form 
$$
E^{a_1}_{k_1} E^{a_2}_{k_2} \ldots E^{a_n}_{k_n} 
\Delta^{b_1}_{k_1} \Delta^{b_2}_{k_2} \ldots \Delta^{b_n}_{k_n} 
\beta(n;k_1,\ldots,k_n)
$$ 
with $a_i, b_i \in \{0,1,2,\ldots\}$.
We can interpret these expressions as sums and differences of cardinalities of certain
subsets of monotone triangles with $n$ rows. For instance, 
$$
\Delta_{k_q} \beta(n;k_1,\ldots,k_n)
$$
is the number of monotone triangles that are strictly increasing in southeast direction and with bottom row $(k_1,\ldots, k_{q}+1,\ldots, k_n)$
such that the $(q-1)$-st part of the $(n-1)$-st row is equal to $k_q$ minus 
the number of monotone triangles that are strictly increasing in southeast direction and 
with bottom row $(k_1,\ldots,k_n)$ such that the $q$-th part of the $(n-1)$-st row is equal to $k_q$. In order to prove 
\eqref{equi}, one has to show that these cardinalities add up to the number of 
monotone triangles. Equivalently, one could follow a similar strategy for the identity
$$
\left( \prod_{1 \le p < q \le n} (\id + E_{k_q} \Delta_{k_p}) \right) 
\alpha(n;k_1,\ldots,k_n) = \beta(n;k_1,\ldots,k_n) 
$$
which is equivalent to Theorem~\ref{maininverse}.

\item This is more a remark than another project: to prove Theorem~\ref{main} I have more or 
less carried out an analysis of the recursion \eqref{rec}.  
I originally started  this analysis when considering a 
somehow reversed situation:
let an $(r,n)$ monotone trapezoid be a monotone triangle with the top $n-r$ rows cut off and 
bottom row $(1,2,\ldots,n)$.
Let $\gamma(r,n;k_1,\ldots,k_{n-r+1})$ denote the number of $(r,n)$ monotone trapzoids with 
prescribed top row $(k_1,\ldots, k_{n-r+1})$. In particular, $\gamma(n,n;k)$ is the number of monotone 
triangles with $n$ rows, bottom row $(1,2,\ldots,n)$ and $k$ as entry in the top row. 
In the bijection between alternating 
sign matrices and monotone triangles, the entry in the top row of the monotone triangle 
corresponds to the column of the unique $1$ in the first row of the alternating sign matrix. 
Thus, $\gamma(n,n;k)$ must be equal to \eqref{refined}.
On the other hand, we can also use \eqref{rec} to compute $\gamma(r,n;k_1,\ldots,k_{n-r+1})$: 
$\gamma(1,n;k_1,\ldots, k_n) = 1$ and 
$$
\gamma(r,n;k_1,\ldots, k_{n-r+1}) = \sum_{(l_1,\ldots, l_{n-r+2})}^{(1,k_1,\ldots, k_{n-r+1},n)} 
\gamma(r-1,n;l_1,\ldots, l_{n-r+2}).
$$
With this extended definition, $\gamma(n,n;k)$ is a polynomial in $k$. In the following we list it 
for $n=1,2,\ldots,6$.
\begin{eqnarray*}
  \gamma(1,1;k)&=&  1 \\
  \gamma(2,2;k)&=&  -1 + 3\,k - k^2 \\
  \gamma(3,3;k)&=&  \frac{1}{12} (48 - 92\,k + 103\,k^2 - 40\,k^3 + 5\,k^4) \\
  \gamma(4,4;k)&=&  \frac{1}{72} (-2160 + 5910\,k - 5407\,k^2 + 2940\,k^3 \\ &&  - 919\,k^4 + 150\,k^5 - 10\,k^6) \\
  \gamma(5,5;k)&=&  \frac{1}{1440} (584640 - 1644072\,k + 1970008\,k^2   \\ && - 1211172\,k^3 + 456863\,k^4 - 111708\,k^5 \\ &&  + 
    17462\,k^6 - 1608\,k^7 + 67\,k^8) \\
 \gamma(6,6;k) &=&  \frac{1}{7560} (-73316880 + 225502200\,k   \\ && - 284097336\,k^2 + 204504097\,k^3 \\&& - 91897169\,k^4 + 27466950\,k^5 \\ && - 5651016\,k^6 + 805518\,k^7 \\&&  - 77646\,k^8 + 4655\,k^9 - 133\,k^{10})
\end{eqnarray*} 
Unfortunately, these polynomials are not equal to \eqref{refined}. (For instance, they do not factor over 
$\mathbb{Z}$.)  They only coincide on the combinatorial range $\{1,2,\ldots, n\}$ 
of $k$. However, it might still be possible to compute $\gamma(n,n;k)$ for general $n$.

Strikingly the degree of $\gamma(n,n;k)$ in $k$ is $2n-2$ as the degree of \eqref{refined}. 
This linear growth is again unexpected because the application of \eqref{rec} can more than double 
a polynomial's degree, see Section~\ref{recursion}. However, one can use Lemma~\ref{grad} and an 
extension of Lemma~\ref{fundASM} to show that, more generally, the degree of 
$\gamma(r,n;k_1,\ldots,k_{n-r+1})$ is $2r-2$ in every $k_i$.

\item Finally we have started to investigate a $q$-version  of the formula in 
Theorem~\ref{main}, i.e. a weighted enumeration of monotone triangles with prescribed 
bottom row $(k_1,\ldots,k_n)$ which reduces to our formula as $q$ tends to $1$.

\end{enumerate}

\end{document}